\documentclass[11pt,a4paper]{amsart}

\usepackage[latin1]{inputenc} \usepackage[T1]{fontenc}

\usepackage{graphicx}

\usepackage{multicol, psfrag, subfigure}

\usepackage{amssymb}
\usepackage{amsthm}

\usepackage[dvipsnames]{xcolor}
\xdefinecolor{forestgreen}{named}{ForestGreen}

\usepackage[all]{xy}
\usepackage{stmaryrd}

\everymath{\displaystyle}

\newtheorem{thm}{Theorem}[section] 
\newtheorem{prop}[thm]{Proposition}
\newtheorem{df}[thm]{Definition} 
\newtheorem{lem}[thm]{Lemma}

\newtheorem{ex}[thm]{Example}
\newtheorem{exs}[thm]{Examples}
\newtheorem{rk}[thm]{Remark}

\newcommand{\R}{\mathbb{R}}

\newcommand{\dd}{\mathrm{d}}

\begin{document}             

\date{March 2018}

\bibliographystyle{plain}

\title[On Legendrian cobordisms and generating functions]{On Legendrian cobordisms and generating functions}

\newcommand{\cs}{$^\dagger$} \newcommand{\cm}{$^\ddagger$}

\author{Ma\"{y}lis Limouzineau}

\newcommand{\nfont}{\fontshape{n}\selectfont}

\address{Mathematishes Institut der Universt\"{a}t zu K\"{o}ln}

\begin{abstract}
This note concerns Legendrian cobordisms in one-jet spaces of functions, in the sense of Arnol'd \cite{Arnold}
 -- consisting of big Legendrian submanifolds between two smaller ones. We are interested in such cobordisms which fit with generating functions, and wonder which structures and obstructions come with this notion. As a central result, we show that the classes of Legendrian concordances with respect to the generating function equipment can be given a group structure. To this construction we add one of a homotopy with respect to generating functions.
\end{abstract}

\maketitle

\section*{Introduction}

The goal of this article is to add the generating function constraint to the Legendrian notion of cobordism initiated by Arnol'd in \cite{Arnold}, and to understand what kind of structure comes out. If the definitions works in the one-jet spaces $J^1 (M)$ for any smooth manifold $M$, our main concern is the study of Legendrian cobordisms between (long) Legendrian knots in the $3$-dimensional standard contact space $J^1(\R)$.  

\begin{figure}[ht]
	\begin{center}
		\includegraphics[scale=0.5]{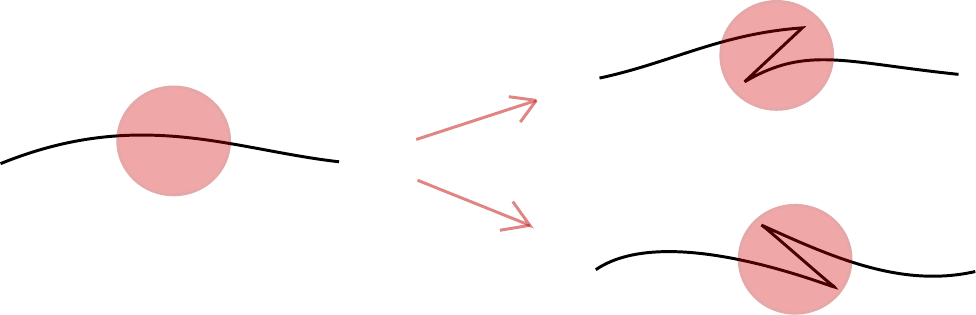}
	\end{center}
	\caption{Stabilizations of a Legendrian knot (front projection).}
		\label{fig1}
\end{figure}

Legendrian knots studied modulo Legendrian isotopy form a more subtle theory than the one of smooth knots. In first approximation, one may combine the smooth theory with the \textit{stabilization operation}, which consists in adding a zig-zag in the front projection (see Figure \ref{fig1}).\\

This operation changes the Legendrian isotopy class of the knot, but not its smooth one. Such a theory can be reduced to the three classical invariants of a Legendrian knot: the topological type, the \textit{Thurston--Bennequin invariant} -- or the self linking number of the knot with a push upward copy of it -- and \textit{Maslov index} -- or he rotation number of its Lagrangian projection. We will say in this introduction that a Legendrian knot which can not be obtained from another one by stabilization is \textit{maximal}. When all Legendrian representatives of a smooth isotopy class of a knot can be listed by doing successive stabilisations to a maximal Legendrian knot, the corresponding topological type is commonly called \textit{simple}.

The Legendrian knot theory does not end here, as there exist different Legendrian isotopy classes of knots with the same classical invariants. It is in general difficult to determine when a Legendrian knot is maximal, and if the corresponding topological type is simple. The first counter-example is the one of Chekanov-Eliashberg \cite{Ch1}, with two distinct Legendrian representatives of the smooth knot type $5_2$, which have the same Thurston-Bennequin invariant and the same Maslov index.

\begin{figure}[ht]
	\begin{center}
		\includegraphics[scale=0.5]{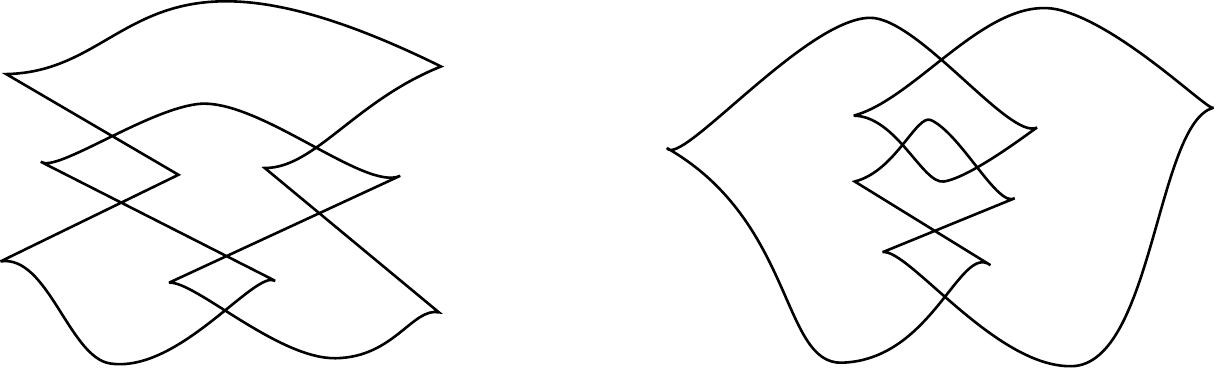}
	\end{center}
	\label{fig2}
	\caption{Wave fronts of the two Chekanov-Eliashberg knots.}
\end{figure}

When working with Legendrian submanifolds in one-jet spaces of functions, generating functions (gf) appear quite automatically as a tool -- here we work more specifically with \textit{generating function} which are \textit{quadratic at infinity} (\textit{gfqi}). Generating functions can be first seen as a Legendrian submanifold factory, constructing many but not all of them. Such constructed Legendrian submanifolds have specific features. In dimension $3$, gf Legendrian knots are maximal (i.e. they have no stabilizations), with Maslov index equal to zero. As a consequence, there exist topological types of knots which do not admit a Legendrian representative equipped with a gf. For instance, one of the two (long) trefoils has no Legendrian representative which has a gfqi, while the other does -- see Proposition \ref{prop1.18}. More generally, Morse theoretic techniques allowed by gf's are used efficiently to provide numerous invariants for Legendrian knots (rulings \cite{PuCh}, \cite{FuRu}, GF-homology \cite{Tr1}, \cite{Tr2}, etc). \\ 

The classical smooth notions of cobordism and concordance have found natural and relevant analogues in contact and symplectic topology. For the last twenty years, the most studied notion has been the one of Lagrangian cobordism between Legendrian submanifolds. To the search for obstructions to the existence of such cobordisms initiated by Baptiste Chantraine in \cite{ChB}, Lisa Traynor and her collaborators successfully added the ingredient of gf's and came out with the notion of \textit{gf-compatible Lagrangian cobordisms} \cite{BST}, \cite{ST}. For instance, Sabloff and Traynor in \cite{ST} found out that there is no gf-compatible Lagrangian cobordisms between the two Chekanov--Eliashberg knots.      

Here we are interested in \textit{Legendrian cobordisms} (see Subsection \ref{ss1.2}), meaning an $(n+1)$-dimensional Legendrian submanifold of $\R^{2n+3}$ between two $n$-dimensional Legendrian submanifolds of $\R^{2n+1}$. 
 	
Because of the extra dimension, this notion is significantly more flexible than the one of (exact) Lagrangian cobordism. The notion of Legendrian cobordism is, in some sense, the extension of Legendrian isotopy study: looking at $1$-parameter families of embedded Legendrian submanifolds, one extends the notion of Legendrian isotopy by adding natural accidents (such as immersion or cobordism moments of Legendrian type, see local descritption Figure \ref{fig5}), whereas the notion of Lagrandian cobordism refines the smooth cobordism theory by adding the symplectic constraint. 
 	
A notable difference is that Lagrangian cobordisms induce embedded smooth cobordisms, while Legendrian cobordisms project on immersed ones.
 	
The gf-equipment extends more naturally to Legendrian cobordisms than to Lagrangian ones. Additionally, it forces the study of $1$-parameter families of Legendrian submanifolds among maximal Legendrian submanifolds. Note that Legendrian knots modulo Legendrian cobordisms have been classified by Arnol'd \cite{Arnold}: two Legendrian knots are Legendrian cobordant if and only if they have the same Maslov index. Our issue is a subclassification, as it concerns Legendrian knots among those with fixed Maslov index equal to zero. 

In this note, we study Legendrian cobordisms with respect to the gfqi-equipment. In particular, we will see that the construction of the smooth concordance group -- the operation here is the connected sum of knots -- has a natural analogue, and this analogue respects the gfqi-equipment.

\begin{thm}{(Theorem \ref{thm2.9}.)} The set of Legendrian knots equipped with a gfqi equivalence class up to Legendrian concordance which is compatible with the gfqi equipment forms a group. 
\end{thm}

The paper is organized as follows: Section 1 is devoted to the necessary tool box of Legendrian geometry. The definitions of Legendrian cobordism (Definition \ref{def1.8}) and gfqi (Definition \ref{def1.13}) are recalled, as well as fundamental results concerning gfqi's (Chekanov Theorem \ref{thm1.16} and Th\'{e}ret--Viterbo Theorem \ref{thm1.15}). We construct the concordance group of Legendrian knots -- \textit{gfqi-concordance} -- in Section 2, and show that the construction fits with gfqi's (Theorem \ref{thm2.9}). Section 3 is devoted to the notion of Legendrian homotopy -- which is a particular case of concordance -- with respect to gfqi's -- \textit{gfqi-homotopy}. We give an additional construction (Proposition \ref{prop3.4}) based on the sum operation \cite{Moi}. \\

The general motives behind those constructions is to discuss the rigidity of the following different notions: Legendrian concordance, Lagrangian concordance, and Legendrian homotopy, with respect to gf's. At the end of this note, some interrogations remains open as far as the author knows : how does the set of gfqi's change along a gfqi-cobordisms? Thanks to the results of Chekanov and Viterbo-Th\'{e}ret, we know for instance that the number of gfqi's does not change along a Legendrian isotopy (see Proposition \ref{prop2.12}), while it can pass from one to infinitely many through a gfqi-cobordism \cite{BST}. What about gfqi-concordance? and gfqi-homotopy?

\section*{Acknowledgements}
I'm grateful to my PhD advisor Emmanuel Ferrand as he supervised my work and gave numerous advices during this note realization. I thank Lisa Traynor, Josh Sabbloff and Samantha Pezzimenti for many interesting discussions concerning Largangian cobordisms and generating functions. I'm also grateful to Sylvain Courtes for enlightening remarks and explanations concerning generating functions. \\

The author is now supported by the grant SFB \textit{Symplectic Structures in Geometry, Algebra and Dynamics} and the University of Cologne, Germany.

\vspace*{1cm}
\section{Preliminaries~: Legendrian things in one-jet spaces}\label{sec1}

\vspace{0.2cm}

\subsection{Contact structure} Consider a smooth manifold $M$ of dimension $m$, endowed with a system of local coordinates $x=(x_1, \dots , x_m)$. We refer to $M$ as the \textit{base space}. The space of $1$-jets of functions based on $M$, $J^1(M)=\R \times T^*M$, is a $(2m+1)$-dimensional manifold. By fixing a Riemannian metric on $M$ -- for $M=\R^m$, we choose the Euclidian metric -- $J^1(M)$ is locally endowed with coordinates $(u,x,y)$, with $u \in \R$ and $y=(y_1, \dots , y_m)$ canonically associated to $x$. It carries a natural contact structure $\xi$, which is the hyperplane field defined as $$\xi = \ker (\dd u -y\dd x)= \ker (\dd u - \sum_{i=1}^m y_i \dd x_i). $$ 

A contact structure is a maximally non integrable distribution: a submanifold everywhere tangent to $\xi$ must have dimension no greater than $m$.

\begin{df}\label{def1}
A \textbf{Legendrian submanifold} $L \subset J^1(M)$ is an $m$-dimensional smooth submanifold of $J^1(M)$ which is everywhere tangent to the contact structure, i.e. $$(\dd u - y \dd x)_{\vert L}\equiv 0 . $$ 
\end{df}

\begin{ex}\label{ex11}
Let $f$ be a smooth function defined on $M$. Its \textbf{$1$-graph}  $$j^1f=\lbrace ( u=f(x),x,y=\partial_x f (x) ) \ \vert \ x \in M\rbrace$$
is an elementary example of a Legendrian embedding of $M$ into $J^1(M)$.
\end{ex}

\begin{rk} The cotangent bundle $T^*M$ is naturally endowed with the standard symplectic form $\omega=\dd( y \dd x)$. Thus, the projection of a Legendrian submanifold $L \subset J^1(M)$ onto $(T^*M,\omega)$ is an immersed exact Lagrangian submanifold.  
\end{rk}

To work with Legendrian geometry, it is convenient to use the \textit{front projection} 
\begin{align*}
\mathbf{p_F} : \ J^1(M) \ & \longrightarrow \R \times M \\
(u,x,y) &\longmapsto \ (u,x) .
\end{align*}

The front projection of a Legendrian submanifold is shortly called its \textbf{front}. We systematically picture Legendrian submanifolds of dimensions $1$ and $2$ via their fronts.\\ 

\begin{exs}${}$\\
$\bullet$ \ If a Legendrian submanifold is the $1$-graph of a function $f$, its front is the graph of $f$.\\
$\bullet$ \ In the case $m=1$, the ambient space has dimension $3$, two types of singularities may appear generically in a front: double points and (right or left) cusps. For $m=2$, a swallow tail may appear (see Figure \ref{fig3}, as well as lines of double points, lines of cusps, \dots (see \cite{2} for an exhaustive and detailed description of two dimensional wave fronts).\\
\end{exs}

\begin{df}
A \textbf{Legendrian isotopy} is a smooth one-parameter family of Legendrian submanifolds.
\end{df}

\begin{ex}${}$ In the case of $m=1$, an analogue of the Reidemeister theorem holds: a Legendrian isotopy of knots or links can be seen as a succession of a finite number of local moves on the wave front. The three types of Legendrian Reidemeister moves are illustrated in Figure \ref{fig5} (a).

\end{ex}

\subsection{Legendrian cobordisms}\label{ss1.2}

The notion of Legendrian cobordism as introduced by Arnol'd in \cite{Arnold} consists of an $(n+1)$-dimensional Legendrian submanifold between two $n$-dimensional Legendrian manifolds. It requires a certain \textit{reduction operation} commonly used in Legendrian (and Lagrangian) geometry.

To obtain this Legendrian notion of cobordism, one may consider $M=N\times [0,1]$, where $N$ is a smooth manifold of dimension $n$, and the $1$-jet space $J^1(N\times [0,1])$ endowed with local coordinates $(u,q,t,p,s)$, where $t \in [0,1]$ and $s$ is the dual coordinate of $t$. It carries the natural contact structure $\xi = \ker (\dd u - p \dd q -t \dd s)$.

\begin{df}
Let $\mathcal{L}$ be a subset in $J^1(N \times [0,1])$, and set $t_0 \in [0,1]$. The \textbf{slice of $\mathcal{L}$ at time $t=t_0$} is the projection of $\mathcal{L} \cap \lbrace t=t_0 \rbrace $ on $J^1(N)$ forgetting $s$. We denote it by $\mathcal{L}_{\lceil t=t_0}$, $$\mathcal{L}_{\lceil t=t_0}= \lbrace (u,q,p) \ \vert \ \exists s \in \R \text{ s.t. } (u,q,t_0,p,s)\in \mathcal{L} \rbrace.$$
\end{df}

As a consequence of Thom's Transverality Lemma \cite{Thom}, one can show that the slice at time $t=t_0$ of a Legendrian submanifold of $J^1(N \times [0,1])$ in generic position is a Legendrian submanifold of $J^1(N)$ for almost every $t_0 \in [0,1]$. At the level of fronts, it appears that a generic front in $\R \times N \times [0,1]$ intersected with the hyperplane $\R \times N \times \lbrace t=t_0 \rbrace$ gives a front in $\R \times N$ (see Figures \ref{fig3} and \ref{fig4}).

\begin{df}\label{def1.8}
A \textbf{Legendrian cobordism} consists of a Legendrian submanifold $\mathcal{L} \subset J^1(N \times [0,1])$ and two Legendrian submanifolds $L_0$ and $L_1 \subset J^1(N)$ such that  $$L_0 = \mathcal{L}_{\lceil t=0} \quad \text{ and } \quad L_1 = \mathcal{L}_{\lceil t=1}  .$$
\end{df}

\begin{ex}
A Legendrian isotopy $(L_t)_{t\in [0,1]}$ can be seen as a Legendrian cobordism $\mathcal{L}$ between $L_0$ and $L_1$ such that each slice $\mathcal{L}_{t=t_0}$, $t_0 \in [0,1]$, is an embedded Legendrian submanifold. 
\end{ex}

\begin{ex}
In the case of Legendrian cobordisms between Legendrian knots, the front projection allows us to make reliable illustrations. The first Reidemeister move for instance can be pictured as slices of a swallow tail (Figure \ref{fig3}). In addition to the three Reidemeister moves, three other local moves can be performed on a wave front to obtain Legendrian cobordisms. One is illustrated in Figure \ref{fig4}.

\begin{figure}[ht]
\begin{center}
\includegraphics[scale=0.5]{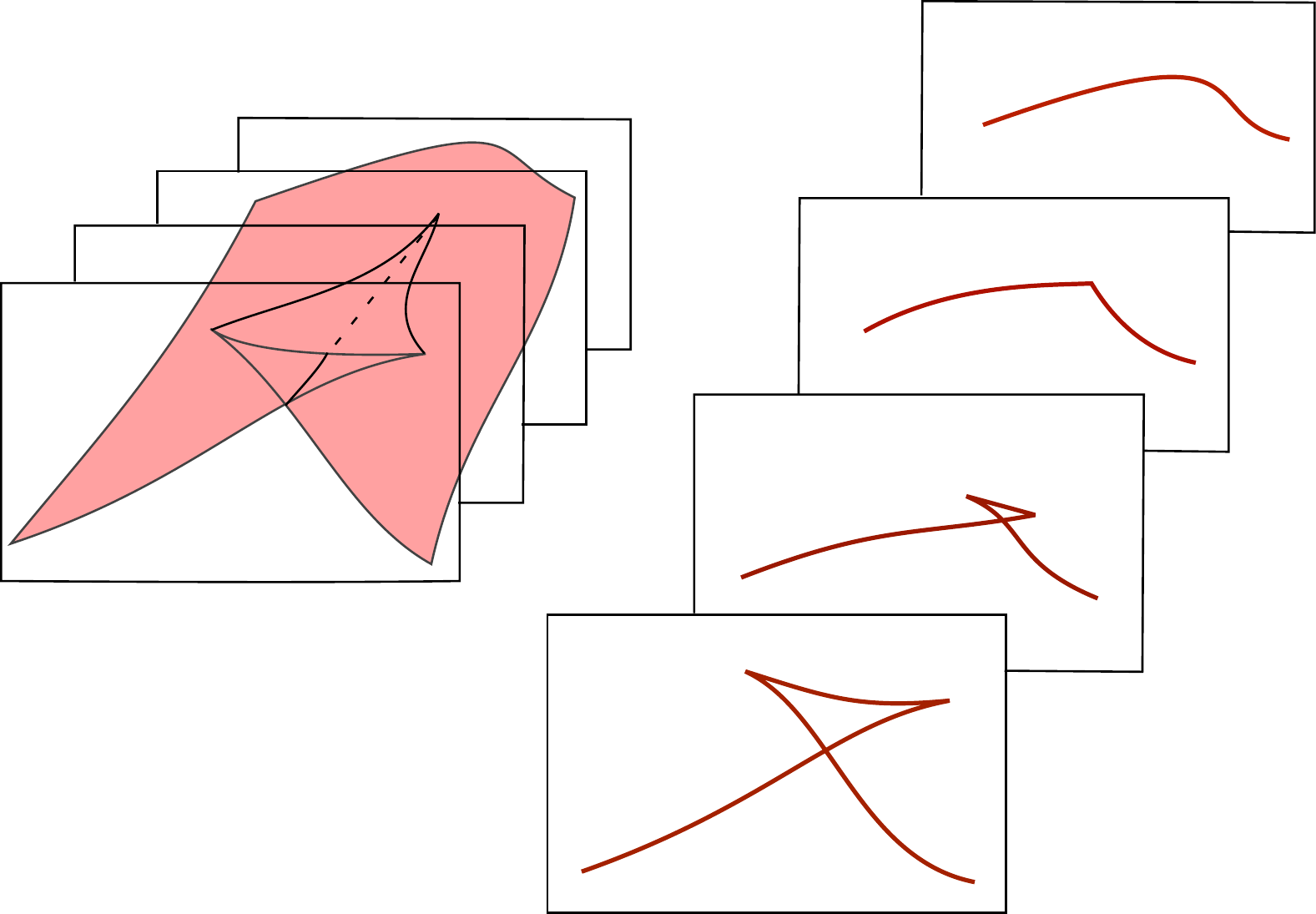}
\end{center}
\caption{The swallow tail and the first Reidemeister move: R I.}
\label{fig3}
\end{figure}

\begin{figure}[ht]
\begin{center}
\includegraphics[scale=0.5]{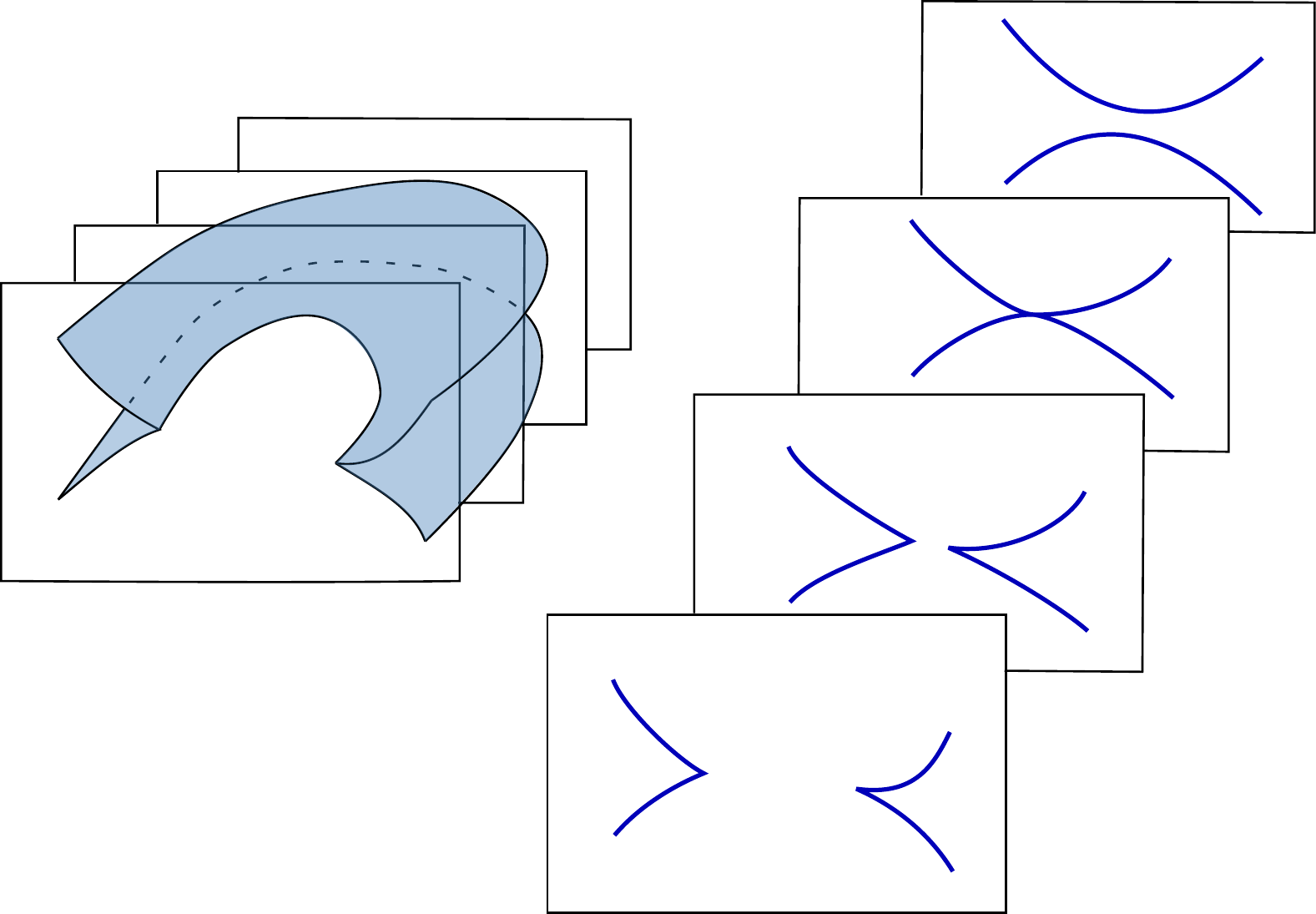}
\end{center}
\caption{Cobordism move: the Legendrian saddle point.}
\label{fig4}
\end{figure}

\begin{figure}[ht]
	\begin{center}
	\subfigure[(a)]{\includegraphics[scale=0.35]{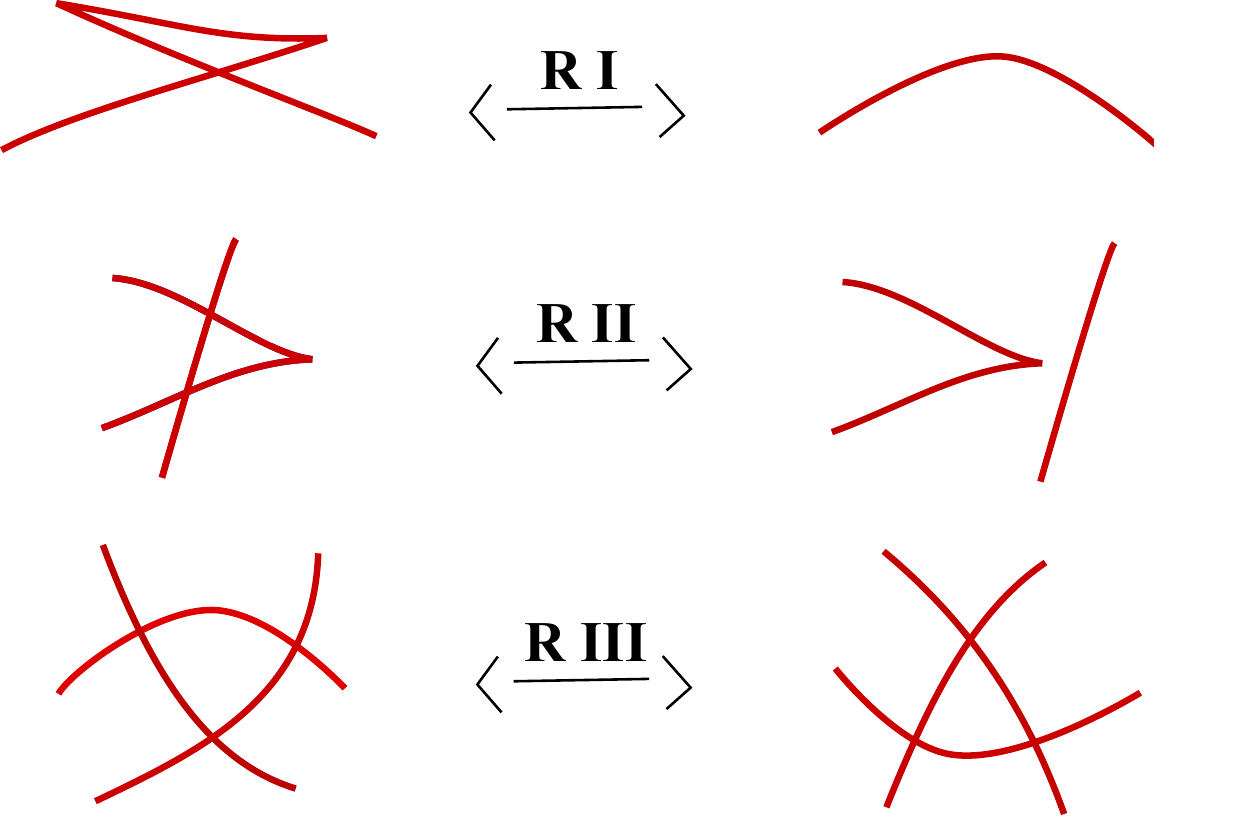}} \qquad \qquad \subfigure[(b)]{\includegraphics[scale=0.35]{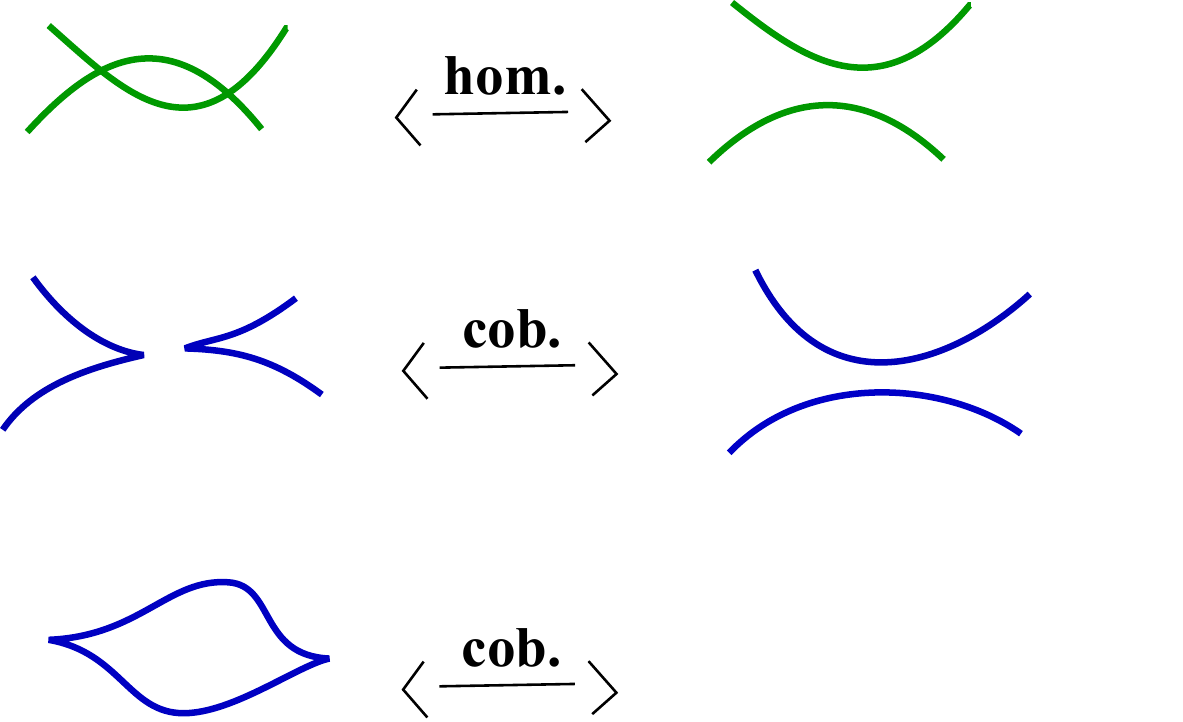}}
	\end{center}
	\caption{The isotopy moves (red), the homotopy move (green) and the cobordisms moves (blue).}
	\label{fig5}
\end{figure}

Comparing Legendrian cobordisms with their smooth analogues, one can observe that the two cobordism moves -- saddle point and spherical -- are completed with a third extra move: the homotopy move (see Figure \ref{fig5}). We will come back to the \textit{homotopy} notion in Section $3$.
\end{ex}

\subsection{Generating functions} The notion of generating function is based on another reduction operation.
Let $k$ be a positive integer. Consider the space $J^1(M\times \R^k)= \R \times T^*(M\times \R^k)$, endowed with coordinates $(u,x,y,w,v)$, where $w \in \R^k$ stands for the extra variable in the base space, and $v \in \R^k$ for the corresponding dual coordinates. Let $\mathcal{L}$ be a subset of $J^1(M\times \R^k)$.

\begin{df}
 The \textbf{contour of $\mathcal{L}$ in the direction of $M$} is the projection of $\mathcal{L} \cap \lbrace v=0  \rbrace$ on $J^1(M)$. 
\end{df} 

As for the slice operation, the contour of a Legendrian submanifold of $J^1(M \times \R^k)$ in generic position is a Legendrian submanifold of $J^1(M)$. In particular, it permits building numerous non-trivial Legendrian submanifolds from $1$-graphs of functions.

\begin{df}
A \textbf{generating function} (gf) for a Legendrian submanifold $L \subset J^1(M \times \R^k)$ is a function $F$ defined on a product $M \times \R^k$ such that $L$ is the contour of the $1$-graph of $F$ in the direction of $M$,
$$L = \lbrace (u,x,y) \  \vert \  \exists w, \partial_w F(x,w)=0, u=F(x,w), y=\partial_x F(x,w) \rbrace .$$
\end{df}

When a Legendrian has a gf, it has infinitely many of them. It is not relevant to distinguish every one of them.
There are two operations which do not change the underlying Morse dynamics of gf: the \textit{stabilization} operation and the \textit{fiberwise diffeomorphism} operation. Together they define a notion of equivalence class for gf such that each invariant constructed from the Morse dynamics of a gf is in fact an invariant of the equivalent class of the gf. 

\noindent $\bullet$ Consider $F$ a gf defined on $M \times \R^k$. Let $k'$ be an integer, and $Q'$ be a non-degenerate quadratic form defined on $\R^{k'}$. Then one may replace $F$ by $F \oplus Q'$, defined on $M \times \R^k \times \R^{k'}$ by  $F \oplus Q (q,w,w')=F(q,w)+Q'(w')$. This operation is called a \textbf{stabilization}.

\noindent $\bullet$ Let $F$ be a gf defined on $M \times \R^k$, and $\Phi$ be a fiberwise diffeomorphism of $M \times \R^k$ -- $\Phi(q,w)=(q,\phi_q(w))$ with $\phi_q$ a diffeomorphism of $\R^k$ for every $q \in M$. One may replace $F$ by the composition $F \circ \Phi$. This operation is called a \textbf{fiberwise diffeomorphism}.  \\

If for two generating functions $F$ and $F'$, there exists a gf $F_0$ such that $F$ and $F'$ descend from $F_0$ by successive stabilizations and fiberwise diffemorphisms operations, then $F$ and $F'$ are declared to be \textbf{equivalent}.

\begin{df}\label{def1.13}
A gf is said \textbf{quadratic at infinity} (\textbf{gfqi}) if it is equivalent to a gf of the form $f+Q$, where $f$ defined on $M\times \R^k$ is a compactly supported gf, and $Q$ is a non-degenerate quadratic form on $\R^k$.
\end{df}

Note that a Legendrian submanifold which is the contour of a gfqi must be equal to the zero section outside a compact set of $J^1(M)$.

\begin{ex} It is possible to realise a Legendrian representative for a long trefoil as the contour of a gfqi $F$ defined on $\R \times \R$. The movie of the one-parameter family of functions $(F(q,.))_{q \in \R}$  is depicted in Figure \ref{figg}. The wave front of the contour of $F$ is also called the \textit{Cerf diagram} of the family $(F(q,.))_{q \in \R}$, \cite{Cerf}.
	
\begin{figure}[ht]
	\begin{center}
		\includegraphics[scale=0.3]{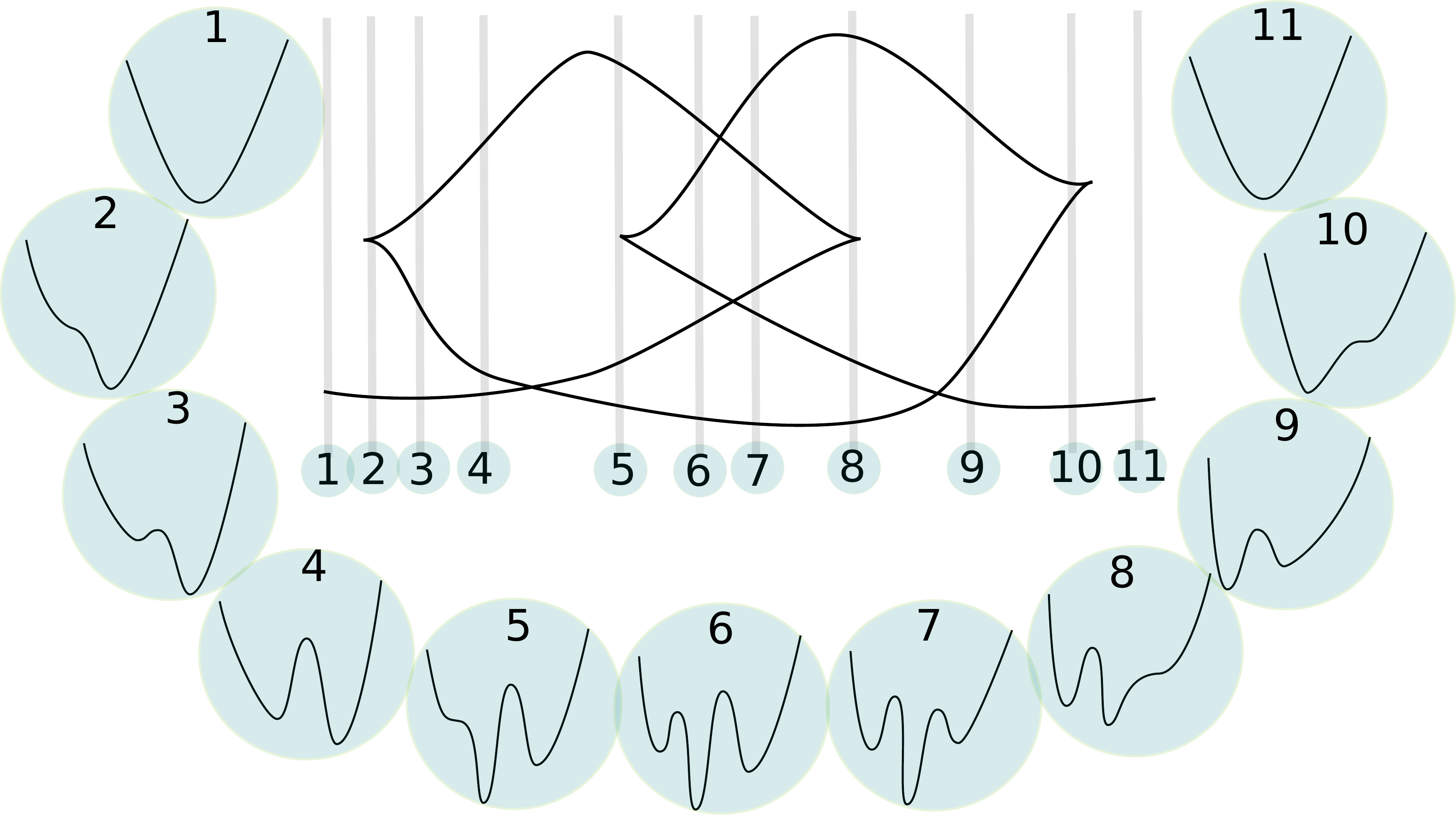}
	\end{center}
	\caption{}
	\label{figg}
\end{figure}
\end{ex}

In Sections 2 and 3, we will use two fundamental results concerning gfqi's. The first one is the persistence of gfqi's under Legendrian isotopies.

\begin{thm}{(Chekanov \cite{Ch}.)}\label{thm1.16}
	
	If $(L_t)_{t\in [0,1]}$ is a Legendrian isotopy, and $L_0$ admits a gfqi $F$, then there exists $\mathcal{F}$ defined on $[0,1]\times M \times \R^k$ such that, for every $t\in [0,1]$, $\mathcal{F}(t, \cdot, \cdot) = F_t$ is a gfqi for $L_t$, and $F_0$ is equivalent to the initial $F$.   
\end{thm}

The second result is the uniqueness of the gfqi class for the zero section.

\begin{thm}{(Th\'{e}ret--Viterbo \cite{Th}, \cite{Vi}.)}\label{thm1.15}
	
If a Legendrian submanifold is isotopic to the zero section, then it admits a unique equivalence class of gfqi's. 
\end{thm}

In particular, the zero section $\mathcal{O}$ admits a unique equivalence class of gfqi's, and we denote it by $F_{\mathcal{O}}$.\\

Back to Legendrian knots, a Morse theoretic argument leads to the observation that a stabilization can never appear when gfqi's -- or any gf's with reasonable behavior at infinity -- are involved. Thus, only maximal Legendrian knots are constructible by gfqi's. Another Morse theoretic observation is that such Legendrian knots must have Maslov index equal to zero. It permits to conclude that all maximal Legendrian knots are not reachable by gfqi's, and show the following result.

\begin{prop}\label{prop1.18}
   The lefthand trefoil does not have any Legendrian long representative that admits a gfqi.
\end{prop}

\begin{figure}[ht]
	\begin{center}
		\includegraphics[scale=0.3]{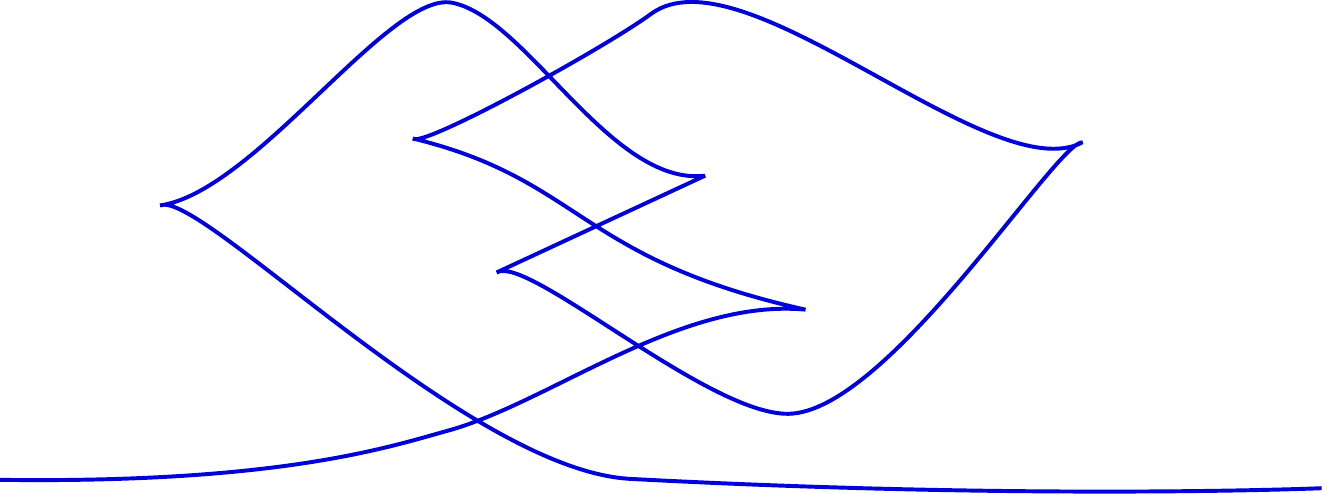}
	\end{center}
	\caption{The maximal (long) Legendrian lefthand trefoil.}
	\label{fig6}
\end{figure}

\noindent \textit{Proof.} It follows from the following fact, proved by Etnyre and Honda in \cite{EtHo}: \textit{the two (topological) trefoil knots are simple.} As a consequence, all long Legendrian representatives of the lefthand trefoil can be obtained by stabilizing the maximal representative depicted in Figure \ref{fig6}, whose Maslov index is not zero. Thus all long Legendrian representatives of the lefthand trefoil with Maslov index zero present a zigzag, therefore can not be realized by a gfqi.
  
\qed

\vspace*{0.5cm}
\section{Concordance of gfqi-knots}\label{sec2}

\subsection{Definitions}

\begin{df}
A \textbf{Legendrian knot} is a connected Legendrian submanifold of $J^1(\R)$.
\end{df}

\begin{df}
A \textbf{gfqi-equipped knot} is a pair $(L,F)$ where $L$ is a Legendrian knot and $F$ is a gfqi equivalence class for $L$.
\end{df}

Thus, a gfqi-equipped knot is equal to the zero section outside a compact set -- it is not a compact but a long knot. Consider the smallest connected open set $U$ in the base space $\R$ such that $L \cap {}^cU= \mathcal{O} \cap {}^cU$. We will refer to $U$ as the \textit{support} of $L$.

\begin{df}
A \textbf{gfqi-cobordism} between two gfqi-equipped knots $(L_0,F_0)$ and $(L_1,F_1)$ consists in a pair $(\mathcal{L},\mathcal{F})$ where $\mathcal{L}$ is a Legendrian cobordism between $L_0$ and $L_1$, and $\mathcal{F}$ is a gfqi equivalence class for $\mathcal{L}$ such that, $$\mathcal{F}_{\vert t=0}=F_0 \quad \text{and} \quad \mathcal{F}_{\vert t=1}=F_1.$$
\end{df}

All gfqi-equipped knots are gfqi-cobordant, as a consequence of the following lemma.

\begin{lem}{(Bourgeois--Sabloff--Traynor \cite{BST}.\footnote{The proof made for generating functions which are \textit{linear at infinity} fits also to the gfqi case.})} \label{lem2.4}

Let $(L,F)$ be a gfqi-equipped knot. There exists a gfqi-cobordism between $(L,F)$ and $(\mathcal{O},F_{\mathcal{O}})$.
\end{lem}

\subsection{The gfqi-concordance group}

\begin{df}
	A \textbf{gfqi-concordance} is a gfqi-cobordism $(\mathcal{L},\mathcal{F})$ such that $\mathcal{L}$ is diffeomorphic to the base space $\R \times [0,1]$.
\end{df}

The notion of gfqi-concordance defines an equivalence relation on the set of gfqi-equipped knots. We denote by $[L,F]$ the equivalence class of a gfqi-equipped knot $(L,F)$ modulo gfqi-concordance.

\subsubsection{The connected sum operation} Let $(L,F)$ and $(L',F')$ be two gfqi-equipped knots.

\begin{figure}[ht]
	\begin{center}
		\includegraphics[scale=0.6]{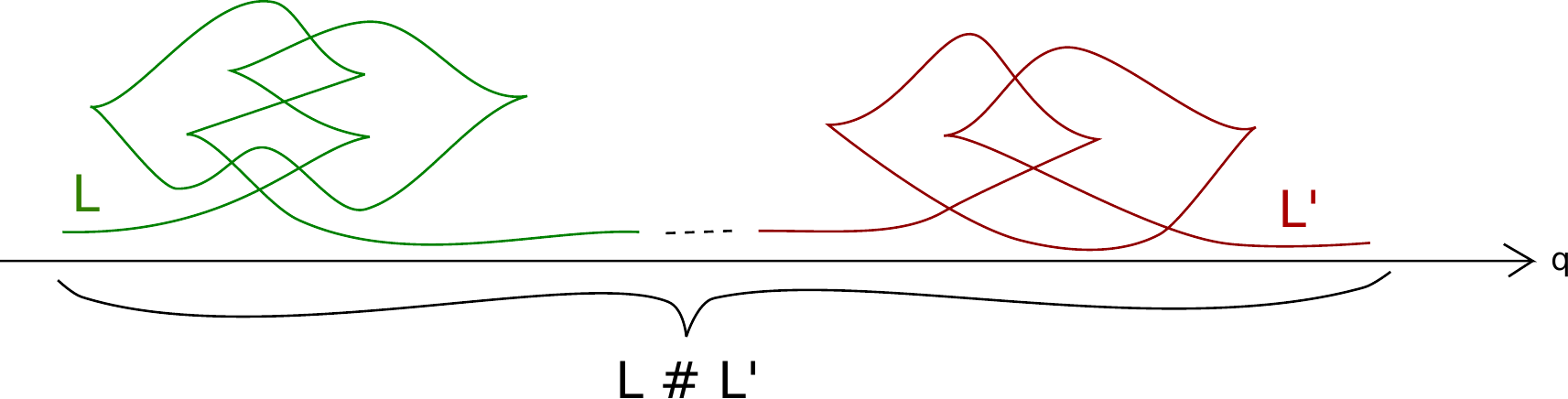}
	\end{center}
	\caption{The connected sum $L\# L'$.}
	\label{fig7}
\end{figure}

\begin{df}\label{def2.6}
	The \textbf{connected sum} of $L$ with $L'$ consists in the concatenation of $L$ then $L'$ (see Figure \ref{fig7}). 
	
\end{df}

The resulting knot, denoted $L \# L'$, has a natural gfqi equipment obtained from $F$ and $F'$:

$$F \# F' \ : \ (q,w,w') \mapsto F(q,w) + F(q,w').$$ 

\begin{rk}
Note that one may have to change $F$ and $F'$ by $F(\cdot +T,\cdot )$ and $F'(\cdot -T, \cdot )$ respectively, with $T \in \R$ large enough, in order to disconnect the supports of $L$ and $L'$ in the base space.
\end{rk}

\begin{rk}
	Reduced to the Legendrian factor, The connected sum operation is commutative modulo isotopy \cite{FuTa}. However, it is not clear that the gfqi-equipment is compatible. Thanks to Chekanov's construction, one can follow this isotopy from $L \#L'$ to $L' \# L$ with a one parameter family of gfqi's, starting with $F \# F'$, but maybe not ending with a gfqi in the same equivalence class as $F' \# F$.
\end{rk}

\begin{df}
	The \textbf{Legendrian mirror} of a Legendrian knot $L$ is the Legendrian knot $\bar{L}$ whose front is the symmetrical of the front of $L$ with respect to the $u$-axis.
\end{df}

If a Legendrian knot $L$ is equipped with a gfqi $F$, we naturally equip its mirror $\bar{L}$ with the gfqi $\bar{F}$ such that: $$\bar{F}(q,\bar{w})=F(-q,\bar{w}).$$

\begin{thm}\label{thm2.9}
Let $(L,F)$ be an gfqi-equipped knot. Then $$[\bar{L},\bar{F}]+[L,F]=[L,F]+[\bar{L},\bar{F}]=[\mathcal{O},F_{\mathcal{O}}].$$
\end{thm}

\noindent \textit{Proof.} We use a well-known construction in Legendrian geometry called \textit{front spinning}. It consists of including the wave front of a Legendrian submanifold in a bigger space, and then rotating it using the additionnal coordinates in order to create a bigger wave front (see Ekholm-Etnyre-Sullivan \cite{EES} and Golovko \cite{G}). 

Here consider $\mathbf{p_F}(L)$ the wave front of $L$, which is one-dimensional. Suppose its support is located in the half space $\lbrace q\geqslant 0 \rbrace$. Let's add the $t$-coordinate and consider $\mathbf{p_F}(L) \cap \lbrace q \geqslant 0 \rbrace$ included in the subspace 

\begin{figure}[ht]
	\begin{center}
		\includegraphics[scale=0.45]{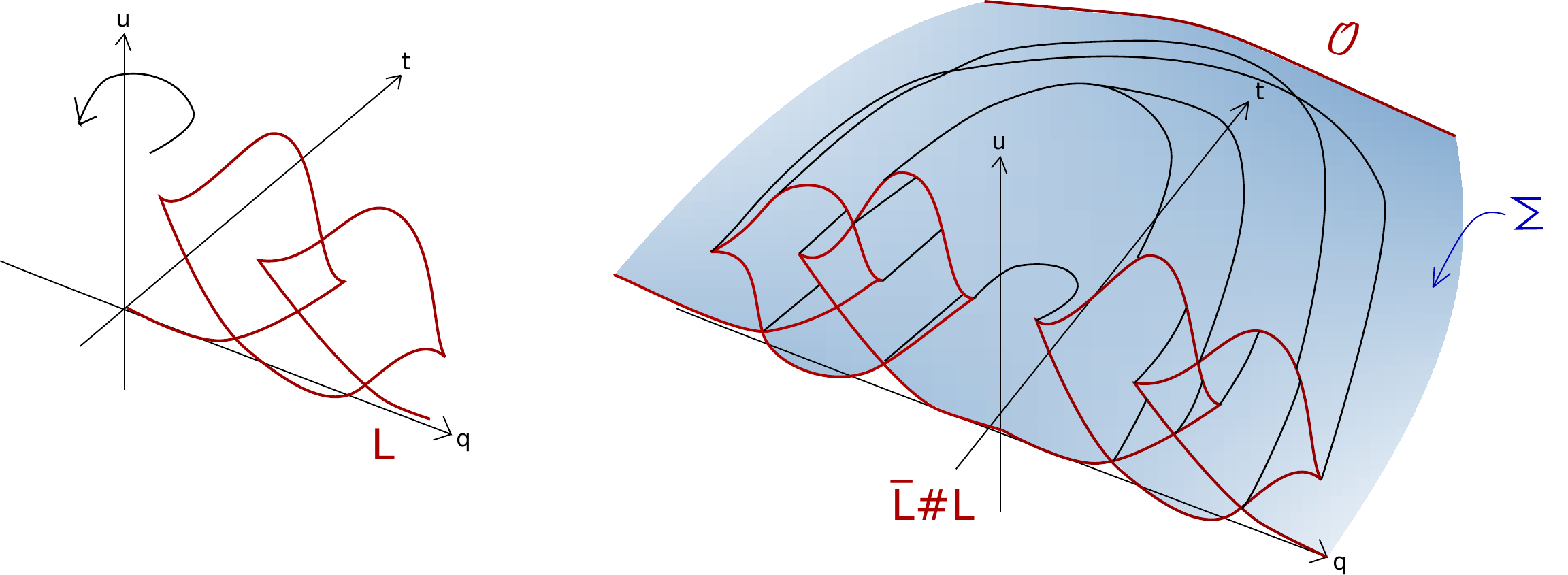}
	\end{center}
	\caption{Construction of the concordance between $\bar{L}\# L$ and $\mathcal{O}$.}
	\label{fig8}
\end{figure}

\noindent $\lbrace (u,q,t) \ \vert \ q \geqslant, t=0 \rbrace \subset \R^2 \times [0,1]$. Then make half a turn around the $u$-axis with this singular curve, and keep the trace of the rotation all along to create a $2$-dimensional object in $\R^2 \times [0,1]$ (see Figure \ref{fig8}).\\

The result of this operation is the wave front of a Legendrian surface $\Sigma L$ in $J^1(\R \times [0,1])$, with a two component boundary. One component lives in the subspace $\lbrace t=0 \rbrace$ and corresponds to the connected sum $\bar{L} \# L$. The second is the subspace $\lbrace (u,q,t,p,s) \ \vert \ u=p=s=0, t=1 \rbrace$, and corresponds to the zero section $\mathcal{O} \subset J^1(\R)$. 

Thus, the Legendrian surface $\Sigma L $ realizes a Legendrian cobordism between $\bar{L} \# L$ and the zero section. As $\Sigma$ has clearly genus zero, it is a concordance.\\

If $F$ is a gfqi for $L$, then $\Sigma L$ has the following gfqi,

$$\mathcal{F}: (q,t,w) \mapsto F(\sqrt{q^2+t^2},w).$$

For $t=1$, the restriction $\mathcal{F}_{\vert t=1} $ is a gfqi for $\mathcal{O}$. Thanks to the unicity result of Viterbo-Th\'{e}ret -- Theorem \ref{thm1.15} -- we obtain for sure the zero section endowed with its unique gfqi-equipment, $(\mathcal{O},\mathcal{F}_{\mathcal{O}})$, at the end of this gfqi-concordance.\footnote{Note that it is not needed to call for Viterbo-Th\'{e}ret theorem here, as this fact follows by construction.}

Replace $F$ by an equivalent form $f+Q$, where $f$ is a gf compactly supported and $Q$ a non-degenerate quadratic form of the $w$-variable. For $t=0$, let us write 
$$\mathcal{F}_{t=0}(q,0,w)=f(-q,w)+f(q,w)+Q(w),$$
and compare with $$\bar{F}\# F (q,\bar{w},w)=f(-q,\bar{w})+Q(\bar{w})+ f(q,w)+Q(w).$$

In \cite{Th}, in order to prove the invariance of uniqueness property under isotopies, Th\'{e}ret was led to prove the following result.

\begin{lem}\label{lem2.10}
	If $(F_t)_{t\in [0,1]}$ is a smooth path of gfqi's which have all the same contour, $L_t=L \ , \ \forall t$ , then $F_0$ and $F_1$ are equivalent. 
\end{lem}   

Thanks to this technical Lemma, it is sufficient to link the equivalent classes of $\mathcal{F}_{t=0}$ and $\bar{F}\# F$ by a path of gfqi's whose contour is constant and equal to $\bar{L}\# L$, in order to conclude. 

Suppose $F$ is defined on $\R\times \R^k$. We replace $\mathcal{F}_{t=0}$ by the equivalent gfqi $F_0$ defined on $\R \times \R^k \times \R^k$ by
$$F_0(q,w,\bar{w})=f(-q,w)+f(q,w)+Q(w)+Q(\bar{w}). $$

Then, let us define the path of gfqi's $(F_t)_{t\in [0,1]}$ by $$F_t(q,w,\bar{w})=f(-q,\cos(\tfrac{\pi}{2} t) w +\sin(\tfrac{\pi}{2} t)\bar{w})+f(q,w) + Q(w)+Q(\bar{w}).$$

One can check that the contour remains constant, and clearly this path links $F_0$ with $F_1= \bar{F} \# F$. \qed\\

\begin{rk}
	We wonder if the gfqi-concordance group is trivial or not. Remind that Chekanov theorem claims that, if $(L_0, F_0)$ is a gfqi-equipped knot, and $(L_t)_{t \in [0,1]}$ is a Legendrian isotopy from $L_0$ to $L_1$, then there exists -- a one parameter family of gfqi's $(F_t)_{t\in [0,1]}$, ending with -- a gfqi $F_1$ for $L_1$ such that $[L_0,F_0]=[L_1,F_1]$. However, the gfqi class $F_1$ for $L_1$ can not be fixed in advance. 
	Moreover, a stronger version of Lemma \ref{lem2.10} holds:
\end{rk}	

\begin{prop}\label{prop2.12}
Let $(F_t)_{t\in [0,1]}$ and $(F_t')_{t\in [0,1]}$ be two smooth paths of gfqi's, such that the corresponding contours $(L_t)_{t\in [0,1]}$ form the same Legendrian isotopy. Suppose $F_0$ and $F_0'$ are equivalent. Then $F_1$ and $F_1'$ are equivalent.
\end{prop}	

\noindent \textit{Proof.} In \cite{Th}, Th\'{e}ret proved that the set of gfqi's forms a Serre fibration over the set of Legendrian submanifolds in $J^1(M)$ which are diffeomorphic to the base space $M$.\footnote{More precisely, it is done for Lagrangian submanifolds in \cite{Th}, and adapted for the Legendrian case in \cite{14}.} Consider the path of gfqi's $F\star F'^{-1}$, formed by composing $(F_t)_{t\in [0,1]}$ with $(F_t')_{t\in [0,1]}$ traveled the other way around. It gives a path from $F_1'$ to $F_1$ passing by $F_0=F_0'$, which projects onto the loop of Legendrian submanifolds formed by composing the path $(L_t)_{t_\in [0,1]}$ with its inverse. It is contractible, and retracts on the constant loop $(L_1)_{t\in [0,1]}$. Thus there exists a retraction by deformation from the path $F \star F'^{-1}$ to a path $(\tilde{F}_t)_{t\in[0,1]}$, with $\tilde{F}_0=F_1'$ and $\tilde{F}_1=F_1$, which satisfies the assumption of Lemma \ref{lem2.10}. \qed\\

	In other words, the number of equivalence classes of gfqi's for a Legendrian knot remains the same along Legendrian isotopies.\\
	
\noindent \textbf{Question:} How does this number change along a gfqi-concordance (or \textit{homotopy}, see next Section)?	
	
\section{A gfqi-homotopy construction}\label{sec21}

This third part is devoted to another cobordism notion with respect to the existence of gfqi's, between Legendrian isotopy and gfqi-concordance. The following construction emphasizes the flexibility that exists among gfqi-concordances.     

\begin{df}\label{def21}
A \textbf{Legendrian homotopy} is a Legendrian cobordism $\mathcal{L}$ such that all the $t_0$-slices $\mathcal{L}_{\lceil t =t_0}, t_0 \in [0,1]$ are in immersed Legendrian submanifolds.     
\end{df}

A Legendrian homotopy is a particular case of a concordance, avoiding cobordisms moments. In terms of local decomposition on a wave front (see Figure \ref{fig5}), a generic Legendrian homotopy can be decomposed using only Reidemester moves I, II, III, and the homotopy move IV.\\

In \cite{EF}, E. Ferrand shows that the homotopy move can be avoided, by replacing it with a sequence of isotopy and cobordism moves. In other words, if there is a Legendrian cobordism between two Legendrian knots, one may change it to another cobordism which projects onto $J^1(\R) \times [0,1]$ on an smooth embedded cobordism.

However, this sequence of isotopy and cobordism moves reveals stabilised knots or links, which are forbidden in our context of working with generating functions. 

\begin{df}
A \textbf{homotopy with a gfqi} is a Legendrian homotopy $\mathcal{L}$ such that there exists a gfqi $\mathcal{F}$ for $\mathcal{L}$.  
\end{df}

Let $(L,F)$ be a gfqi-equipped knot. Remind that Theorem \ref{thm2.9} implies that there exists a concordance between the connected sum $\bar{L} \# L$ and the zero-section $\mathcal{O}$ which has a gfqi. 

\begin{prop}\label{prop3.4}
	Let $L$ be a Legendrian knot having a gfqi $F$. Then there exists a homotopy with a gfqi between $L \# \bar{L} \# L$ and $L$, $$L \# \bar{L} \# L \underset{\text{hom. gfqi}}{\sim} L.$$  
\end{prop}

\noindent \textit{Proof.} The construction is based on the \textit{sum} operation, defined in \cite{Moi}. Starting with two generating functions $F_1$ and $F_2$ for respectively two Legendrian submanifolds $L_1$ and $L_2$, one obtains a third Legendrian object, denoted $L_1 \underset{\smile}{+} L_2$, by summing up the gf's over the base space. It is generically an embedded Legendrian submanifold, which admits a gf $F_1 \underset{\smile}{+} F_2$ defined as $$F_1 \underset{\smile}{+} F_2(q, w_1,w_2)=F_1(q,w_1)+F_2(q,w_2).$$

If $L_1$ and $L_2$ have disjoint supports over the base space, then the sum operation is nothing else than the connected sum operation, Definition \ref{def2.6}. Note also that, if $F_1$ and $F_2$ are quadratic at infinity, so is $F_1 \underset{\smile}{+} F_2$.

\begin{figure}[ht]
	\begin{center}
		\includegraphics[scale=0.9]{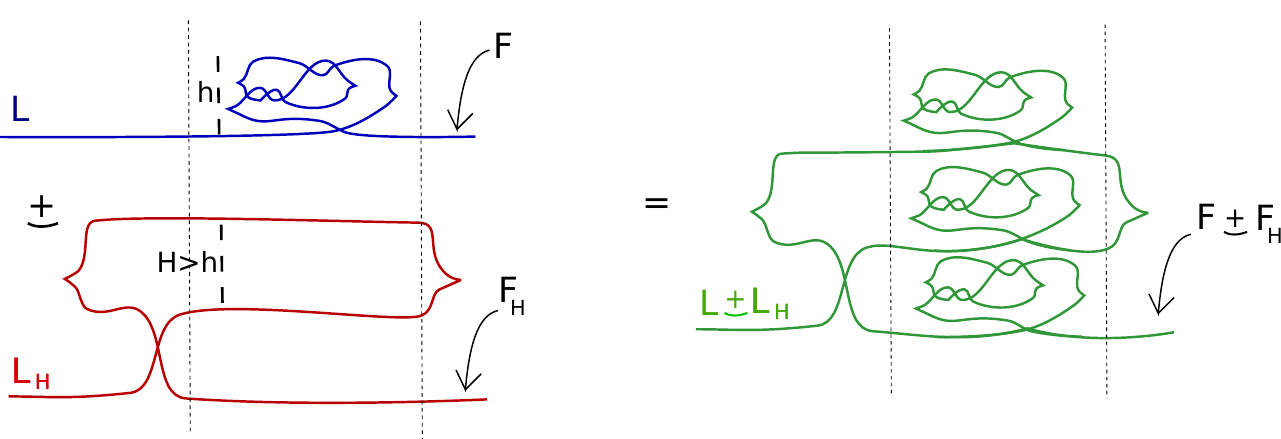}
	\end{center}
	\caption{}
	\label{fig9}
\end{figure}

Consider now the Legendrian submanifold $L_H$, isotopic to the zero section, described by its wave front in red in Figure \ref{fig9}. It must lie in such a way that over the support of $L$, $L_H$ consists of three horizontal strings.  We will adjust the spacing between the two upper horizontal strings. If the knot $L$ has is height equal to $h$, then we start with a spacing $H$ such that $H>h$. Summing up $L$ with $L_H$ then gives three copies of $L$ attached at each level of $L_H$. The resulting Legendrian submanifold is isotopic to the connected sum $L \# \bar{L} \# L$.

\begin{figure}[ht]
	\begin{center}
		\includegraphics[scale=0.9]{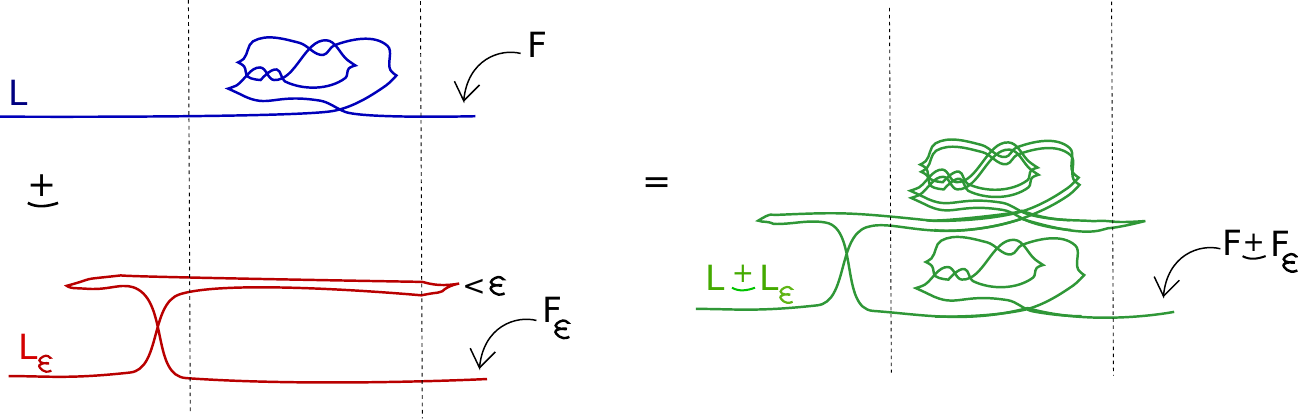}
	\end{center}
	\caption{}
	\label{fig10}
\end{figure}

We then decrease the spacing between the two upper strings of $L_H$ from $H$ to $\epsilon$, for a sufficiantly small enough $\epsilon$, see Figure \ref{fig10}. The Legendrian submanifold $L \underset{\smile}{+} L_{\epsilon}$ at the end of the procedure is isotopic to $L$. 

The $1$-parameter family $(L \underset{\smile}{+} L_{t})_{t\in [\epsilon, H]}$ is a homotopy, which at the end gives us a homotopy from $L \# \bar{L} \# L$ to $L$.

There is a simple $1$-parameter family of gfqi's $(F_t)_{t \in [\epsilon, H]}$ for $(L_t)_{t\in [\epsilon , H]}$, which can be extended in a gfqi for the whole homotopy from $L \# \bar{L} \# L$ to $L$ thanks to Chekanov's theorem.

\qed\\

The different constructions of Legendrian cobordisms with gfqi's should be compared with the notions of gf-compatible Lagrangian cobordisms of L. Traynor and collaborators (J. Sabloff, S. Pezzimenti). Note that the homotopy move corresponds to the immersed points of gf-compatible Lagrangian cobordisms. As far as we know, it is clear that homotopy moves can not be avoided in all gf-compatible (immersed) Lagrangian cobordisms. For instance, there is no embedded gf-compatible Lagrangian cobordism between the two Chekanov's knots, \cite{ST}. 

On the one hand, one wonders when homotopy moves can be removed with respect to gf's to create embedded Lagrangian cobordisms. 

On the other hand, we ask if it is possible to systematically remove cobordism moments rather than immersed points. In other words, using the terminology of this note, is there any obstruction for the existence of a gfqi-homotopy between two gfqi-equipped knots (i.e. with respect to the gfqi-equipment)? Note that a positive answer would imply that the gfqi-concordance group is trivial.\\

\noindent \textbf{Question:} Let $L$ be a gfqi-knot such that there exist two different gfqi classes $F_1$ and $F_2$ whose contour is $L$. Does it exist a gfqi-isotopy from $(L,F_1)$ to $(L,F_2)$? a gfqi-concordance? a gfqi-homotopy?

\bibliography{biblio}

\vspace*{2cm}
\end{document}